\documentclass[12pt,leqno]{article}
\usepackage{latexsym}
\usepackage{bbm}
\usepackage{amssymb,amsmath}
\begin{titlepage}
\title{On the large sieve with square moduli}
\author{Stephan Baier}
\date{22.08.05}
\end{titlepage}
\begin{document}
\maketitle
$ $\\
$ $\\
{\bf Address of the author:}\medskip\\ 
Stephan Baier\\ 
Jeffery Hall\\ 
Department of Mathematics and Statistics\\ 
Queen`s University\\
University Ave\\
Kingston, Ontario, Canada\\
K7L 3N6\medskip\\ 
e-mail: sbaier@mast.queensu.ca
\newpage 
$ $\\
{\bf Abstract:} We prove an estimate for the large sieve with square moduli which improves a recent result of L. Zhao. Our method uses an idea of D. Wolke and some results from Fourier analysis.
 
\section{Introduction}
Throughout this paper, we reserve the symbols 
$c_i$ $(i=1,2,...)$ for absolute positive
constants, and the symbol $\varepsilon$ for an arbitrary (small) positive 
constant. 
Further, we suppose that 
$(a_n)$ is a sequence of complex numbers and that $Q,N\ge 1$. We set
$$
S(\alpha):=\sum\limits_{n\le N} a_n e(n\alpha)
$$
and 
$$
Z:=\int\limits_0^1 \vert S(\alpha) \vert^2 {\rm d} \alpha =
\sum\limits_{n\le N} \vert a_n\vert^2.
$$

In its modern form, the large sieve is an inequality connecting 
a discrete and the continuous
mean value $Z$ of the trigometrical polynomial $S(\alpha)$, {\it i.e.} an 
inequality of the form 
$$
\sum\limits_{r=1}^R \left\vert S\left(\alpha_r\right)\right\vert^2
\le f(N;\alpha_1,...,\alpha_r)Z.
$$
One formulation of the large sieve is as follows.\\ 

{\bf Theorem 1:} \begin{it} Let $\left(\alpha_r\right)_{r\in\mathbbm{N}}$ be a
sequence of real numbers. Suppose that $0<\Delta\le 1/2$ and 
$R\in \mathbbm{N}$. Put 
\begin{equation}
K(\Delta):=\max\limits_{\alpha\in \mathbbm{R}} 
\sum\limits_{\scriptsize \begin{array}{cccc} r=1\\
\vert\vert \alpha_r -\alpha\vert\vert\le \Delta \end{array}}^R 1,\label{1}
\end{equation}
where $\vert\vert x \vert\vert$ denotes the distance of a real $x$
to its closest integer.
Then 
$$
\sum\limits_{r=1}^R \left\vert S\left(\alpha_r\right)\right\vert^2
\le c_{1}K(\Delta)(N+\Delta^{-1})Z.
$$
\end{it}

The above Theorem 1 is an immediate consequence of Theorem 2.11 in 
\cite{Lem}.

In many applications, the sequence $\alpha_1,...,\alpha_R$ consists of  
Farey fractions. If $\alpha_1,...,\alpha_R$ is the sequence of all
fractions $a/q$ with $1\le a\le q$, $(a,q)=1$ and $q\le Q$, then the 
above Theorem 1 implies, on choosing $\Delta:=1/Q^2$, that 
$$
\sum\limits_{q\le Q} \sum\limits_{\scriptsize \begin{array}{cccc} 
a=1\\(a,q)=1\end{array}}^q \left\vert S\left(\frac{a}{q}\right)\right\vert^2
\ll (N+Q^2)Z.
$$
This is the classical large sieve inequality of Bombieri \cite{Cla}. 

Recently, L. Zhao \cite{Zha} considered the case when the moduli
$q$ are squares. A careful investigation of the term $K(\Delta)$ for 
this situation led him to the estimate
\begin{equation}
\sum\limits_{q\le Q} \sum\limits_{\scriptsize \begin{array}{cccc} 
a=1\\(a,q)=1\end{array}}^{q^2}
\left\vert S\left(\frac{a}{q^2}\right)\right\vert^2
\ll (\log 2Q)\left(Q^3+(N\sqrt{Q}+\sqrt{N}Q^2)N^{\varepsilon}\right)Z. 
\label{3}
\end{equation}
In \cite{Bai} we proved that the middle term $N\sqrt{Q}$ on the right-hand
side of (\ref{3}) can be replaced by $N$, which gives an improvement of
(\ref{3}) if $Q\ll N^{1/3-\varepsilon}$.

In the present paper we prove\\

{\bf Theorem 2:} \begin{it} We have
\begin{equation}
\sum\limits_{q\le Q} \sum\limits_{\scriptsize \begin{array}{cccc} 
a=1\\(a,q)=1\end{array}}^{q^2}
\left\vert S\left(\frac{a}{q^2}\right)\right\vert^2
\ll (\log Q)N^{\varepsilon}\left(Q^3+N^{5/4}\right)Z. 
\label{4}
\end{equation}\end{it} 

This bound is sharper than (\ref{3}) if 
$N^{3/8+\varepsilon}\ll Q \ll N^{1/2-\varepsilon}$.

To establish Theorem 2, we combine a method of D. Wolke \cite{Wol} with
some standard 
tools from harmonic analysis,
like the Poisson summation formula and 
bounds for exponential integrals. We also use a bound for  
quadratic Gau\ss{} sums. 
    
\section{Counting Farey fractions in short intervals}
To prove Theorem 2, we shall use the general large sieve bound given in 
Theorem 1.

In the sequel, suppose that $Q_0\ge 1$, and let  
$\alpha_1,...,\alpha_R$ be the sequence
of Farey fractions $a/q^2$ with $Q_0\le q^2\le 2Q_0$, 
$1\le a\le q^2$ and $(a,q)=1$. 
Suppose that $\alpha\in \mathbbm{R}$ and $0<\Delta\le 1/2$. 
Put
$$
I(\alpha):=[\alpha-\Delta,\alpha+\Delta]
\ \ \mbox{ and }\ \ 
P(\alpha):= \sum\limits_{\scriptsize \begin{array}{cccc} 
Q_0\le q^2\le 2Q_0\smallskip\\ (a,q)=1\smallskip\\
a/q^2\in I(\alpha) \end{array}} 1.
$$
Then we have
$$
K(\Delta)=\max\limits_{\alpha\in \mathbbm{R}} P(\alpha),\label{60}
$$
where $K(\Delta)$ is defined as in (\ref{1}). Therefore, the proof of Theorem 
2 reduces to estimating $P(\alpha)$ for all 
$\alpha\in \mathbbm{R}$ and choosing the parameter $\Delta$ 
appropriately.   

To estimate $P(\alpha)$, we begin with a method of D. Wolke \cite{Wol}. Let
\begin{equation}
\tau:=\frac{1}{\sqrt{\Delta}}.\label{P1}
\end{equation}
Then, by Dirichlet's approximation theorem, $\alpha$ can be written in the form
\begin{equation}
\alpha=\frac{b}{r}+z, \ \ \mbox{ where }\ \  r\le \tau,\ (b,r)=1,\ 
\vert z\vert \le \frac{1}{r\tau}.\label{P2}
\end{equation}
Thus, it suffices to estimate $P(b/r+z)$ for all $b,r,z$ satisfying
(\ref{P2}).

We further note that we can restrict ourselves to the case when 
\begin{equation}
z\ge \Delta.\label{P4}
\end{equation}
If $\vert z\vert<\Delta$, then
$$ 
P(\alpha)\le P\left(\frac{b}{r}-\Delta\right)+P\left(\frac{b}{r}+\Delta\right).
$$ 
Furthermore, we have
$$
\Delta=\frac{1}{\tau^2}\le \frac{1}{r\tau}.
$$ 
Therefore
this case can
be reduced to the case $\vert z\vert=\Delta$. 
Moreover, as $P(\alpha)=P(-\alpha)$, we can choose
$z$ positive. So we can assume (\ref{P4}).   

Summarizing the above observations, we deduce\\

{\bf Lemma 1:} \begin{it} We have  
\begin{equation}
K(\Delta)\le 2\max\limits_{\scriptsize \begin{array}{cccc} 
r\in \mathbbm{N} \\ r\le 1/\sqrt{\Delta} \end{array}} 
\max\limits_{\scriptsize \begin{array}{cccc} b\in \mathbbm{Z}\\ 
(b,r)=1\end{array}}
\max\limits_{\Delta\le z\le \sqrt{\Delta}/r} P\left(\frac{b}{r}+z\right).
\end{equation} \end{it}

The next lemma provides a first estimate for $P\left(b/r+z\right)$.\\
       
{\bf Lemma 2:} \begin{it}
Suppose that the conditions (\ref{P1}), (\ref{P2}) and (\ref{P4}) are 
satisfied. Suppose further that 
\begin{equation}
\frac{Q_0\Delta}{z}\le \delta\le Q_0. \label{P10}
\end{equation}
Then,
\begin{eqnarray}
\label{PP}\\
& & P\left(\frac{b}{r}+z\right)\nonumber\\ &\le& c_1\left(1+
\frac{1}{\delta} \int\limits_{Q_0}^{2Q_0} \left(
\sum\limits_{\sqrt{y}-c_2\delta/\sqrt{Q_0}\le 
q \le \sqrt{y}+c_2\delta/\sqrt{Q_0}}
\sum\limits_{\scriptsize 
\begin{array}{cccc} (y-4\delta)rz\le m \le (y+4\delta)rz\smallskip
\\ m \equiv -bq^2\mbox{ mod } r\smallskip\\ m\not=0\end{array}} 1 \right)
{\rm d}y\right).\nonumber
\end{eqnarray}
\end{it}\\

{\bf Proof:} By $\delta\le Q_0$, we have 
\begin{equation}
P(\alpha)\le \frac{1}{\delta} \int\limits_{Q_0}^{2Q_0} P(\alpha,y,\delta) 
\ {\rm d} y,\label{P6}
\end{equation}
where  
$$
P(\alpha,y,\delta):=\sum\limits_{\scriptsize \begin{array}{cccc} 
y-\delta\le q^2\le y+\delta\smallskip\\ (a,q)=1\smallskip\\
a/q^2\in I(\alpha) \end{array}} 1.
$$
Now, for  
$$
y-\delta\le q^2\le y+\delta,\ \ \ \ \ \ (a,q)=1,\ \ \ \ \ \ 
\frac{a}{q^2}\in I(\alpha),
$$
we have $q^2(\alpha-\Delta)\le a\le q^2(\alpha+\Delta)$ or, by (\ref{P2}) and 
(\ref{P4}), $(y-\delta)r(z-\Delta)\le ar-bq^2\le (y+\delta)r(z+\Delta)$. 
If $ar-bq^2=0$, then 
$r=q^2$ since $(a,q)=1=(b,r)$.
Hence,
\begin{equation}
P(\alpha,y,\delta)\le \nu(y)+
\sum\limits_{\sqrt{y-\delta}\le q\le \sqrt{y+\delta}}
\sum\limits_{\scriptsize 
\begin{array}{cccc} (y-\delta)r(z-\Delta)\le m \le (y+\delta)r(z+\Delta)
\smallskip\\ m \equiv -bq^2\mbox{ mod } r\smallskip\\ m\not=0
\end{array}} 1, \label{P11}
\end{equation}
where 
$$
\nu(y):=\left\{\begin{array}{llll} 1, & \mbox{ if } y-\delta\le r\le y+\delta,
\\ \\
0, & \mbox{ otherwise.}\end{array}\right.
$$
Whenever $1\le Q_0\le y\le 2Q_0$ and $\delta\le Q_0$,  
we have, by Taylor's formula, $\sqrt{y}-c_2\delta/\sqrt{Q_0}\le 
\sqrt{y-\delta}\le \sqrt{y+\delta} \le \sqrt{y}+c_2\delta/\sqrt{Q_0}$ 
for a suitable positive constant $c_2$. Furthermore, by (\ref{P10}), we 
have $(y-4\delta)rz\le 
(y-\delta)r(z-\Delta)\le (y+\delta)r(z+\Delta)\le (y+4\delta)rz$. Thus,
(\ref{P11}) implies 
\begin{equation}
P(\alpha,y,\delta)\le \nu(y)+\sum\limits_{\sqrt{y}-c_2\delta/\sqrt{Q_0}\le 
q \le \sqrt{y}+c_2\delta/\sqrt{Q_0}}
\sum\limits_{\scriptsize 
\begin{array}{cccc} (y-4\delta)rz\le m \le (y+4\delta)rz
\smallskip\\ m \equiv -bq^2\mbox{ mod } r\smallskip\\ m\not=0\end{array}} 1. 
\label{P12}
\end{equation}
Combining (\ref{P6}) and (\ref{P12}),
we obtain (\ref{PP}). $\Box$

\section{Estimation of $P(b/r+z)$ - first way}
In this section we use some tools from harmonic analysis to establish
the following bound.\\ 

{\bf Theorem 3:} \begin{it} Suppose that the conditions 
(\ref{P1}), (\ref{P2}) and (\ref{P4}) are satisfied.  
Then,

\begin{equation}
P\left(\frac{b}{r}+z\right)\le 
c_4\Delta^{-\varepsilon}\left(Q_0^{3/2}\Delta
+Q_0^{1/2}\Delta r^{-1/2}z^{-1}+\Delta^{-1/4}\right).\label{QQ}
\end{equation}\end{it}

To derive Theorem 3 from Lemma 2, we need the following standard results from 
Fourier analysis.\\ 

{\bf Lemma 3:} (Poisson summation formula, \cite{Bum}) \begin{it}  
 Let $f(X)$ be a complex-valued 
function on the real numbers that is piecewise continuous with only finitely 
many discontinuities and for all real numbers $a$ satisfies
$$
f(a)=\frac{1}{2}\left(\lim\limits_{x\rightarrow a^-} f(x) +
\lim\limits_{x\rightarrow a^+} f(x)\right).
$$
Moreover, suppose that $f(x)\ll c_5(1+\vert x\vert)^{-c}$ for some $c>1$.
Then,
$$
\sum\limits_{n\in \mathbbm{Z}} f(n) = \sum\limits_{n\in \mathbbm{Z}} \hat{f}(n),
$$
where 
$$
\hat{f}(x):=\int\limits_{-\infty}^{\infty} f(y)e(xy) {\rm d}y,
$$
the Fourier transform of $f(x)$. \end{it}\\

{\bf Lemma 4:} (see \cite{Zha}, for example) \begin{it}
For $x\in \mathbbm{R}\setminus \{0\}$ define
$$
\phi(x):=\left(\frac{\sin \pi x}{2x}\right)^2.
$$
Set 
$$
\phi(0):=\lim\limits_{x\rightarrow 0}\phi(x)=\frac{\pi^2}{4}.
$$ 
Then $\phi(x)\ge 1$ for $\vert x\vert \le 1/2$, and
the Fourier transform of the function $\phi(x)$ is
$$
\hat{\phi}(s)=\frac{\pi^2}{4}\max\{1-\vert s\vert, 0\}.
$$\end{it}

{\bf Lemma 5:} (see Lemma 3.1. in \cite{Kol}) \begin{it} Let $F$ $:$ 
$[a,b]\rightarrow \mathbbm{R}$ be twice differentiable. 
Assume that $\vert F^{\prime}(x)\vert \ge u>0$ 
for all $x\in [a,b]$. Then,
$$
\left\vert \int\limits_{a}^{b} e^{iF(x)} {\rm d}x\right\vert \ll
\frac{c_6}{u}.
$$
\end{it}

{\bf Lemma 6:} (see Lemma 4.3.1. in \cite{Bru}) \begin{it} Let $F$ $:$ 
$[a,b]\rightarrow \mathbbm{R}$ be twice continuously differentiable. 
Assume that $\vert F^{\prime\prime}(x)\vert \ge u>0$ 
for all $x\in [a,b]$. Then,
$$
\left\vert \int\limits_{a}^{b} e^{iF(x)} {\rm d}x\right\vert \le
\frac{c_7}{\sqrt{u}}.
$$\end{it}

We shall also need the following estimate for quadratic Gau\ss{} sums.\\

{\bf Lemma 7:} (see page 93 in \cite{Kol}) \begin{it}
Let $c\in \mathbbm{N}$, $k,l\in 
\mathbbm{Z}$ with $(k,c)=1$. Then,
$$
\sum\limits_{d=1}^r e\left(\frac{kd^2+ld}{c}\right) \le \sqrt{2c}.
$$ 
\end{it}\\

{\bf Proof of Theorem 3:} By Lemma 4, 
the double sum on the right-hand side of (\ref{PP}) can be estimated by
\begin{eqnarray}
& & \sum\limits_{\sqrt{y}-c_2\delta/\sqrt{Q_0}\le 
q \le \sqrt{y}+c_2\delta/\sqrt{Q_0}}
\sum\limits_{\scriptsize 
\begin{array}{cccc} (y-4\delta)rz\le m \le (y+4\delta)rz
\smallskip\\ m \equiv -bq^2\mbox{ mod } r\smallskip
\\ m\not=0\end{array}} 1 \label{P13}\\
&\le& \sum\limits_{q\in \mathbbm{Z}} \ 
\phi\left(\frac{q-\sqrt{y}}{2c_2\delta/\sqrt{Q_0}} 
\right) \sum\limits_{\scriptsize \begin{array}{cccc} 
m\in \mathbbm{Z}\smallskip\\ m\equiv -bq^2 \mbox{ mod } r\end{array}} 
\phi\left(\frac{m-yrz}{8\delta rz}\right) {\rm d}y. \nonumber
\end{eqnarray}  

Using Lemma 3 after a linear change of variables, we transform the inner sum 
on the right-hand side of (\ref{P13}) into 
$$
\sum\limits_{\scriptsize \begin{array}{cccc} 
m\in Z\smallskip\\ m\equiv -bq^2 \mbox{ mod } r\end{array}} 
\phi\left(\frac{m-yrz}{8\delta rz}\right) = 8\delta z
\sum\limits_{j\in \mathbbm{Z}} e\left(\frac{jbq^2}{r}+jyz\right)
\hat{\phi}(8j\delta z).
$$ 
Therefore, we get for the double sum on the right-hand side of (\ref{P13})
\begin{eqnarray}
& & \sum\limits_{q\in \mathbbm{Z}} \ 
\phi\left(\frac{q-\sqrt{y}}{2c_2\delta/\sqrt{Q_0}} 
\right) \sum\limits_{\scriptsize \begin{array}{cccc} 
m\in \mathbbm{Z}\smallskip\\ m\equiv -bq^2 \mbox{ mod } r\end{array}} 
\phi\left(\frac{m-yrz}{8\delta rz}\right)\label{F1} \\
&=& 8\delta z \sum\limits_{j\in \mathbbm{Z}} e(jyz)\hat{\phi}(8j\delta z) 
\sum\limits_{d=1}^{r^*} e\left(\frac{j^*bd^2}{r^*}\right)
\sum\limits_{\scriptsize \begin{array}{cccc} 
k\in \mathbbm{Z}\\ k\equiv d \mbox{ mod } r^*\end{array}}
\phi\left(\frac{k-\sqrt{y}}{2c_2\delta/\sqrt{Q_0}}\right),\nonumber 
\end{eqnarray}
where $r^*:=r/(r,j)$ and $j^*:=j/(r,j)$. Again using Lemma 3 after a linear 
change of variables, we transform the inner sum 
on the right-hand side of (\ref{F1}) into

\begin{equation}
\sum\limits_{\scriptsize \begin{array}{cccc} 
k\in \mathbbm{Z}\\ k\equiv d \mbox{ mod } r^*\end{array}}
\phi\left(\frac{k-\sqrt{y}}{2c_2\delta/\sqrt{Q_0}}\right)=
\frac{2c_2\delta}{r^*\sqrt{Q_0}} \sum\limits_{l\in \mathbbm{Z}} 
e\left(l\cdot\frac{d-\sqrt{y}}{r^*}\right)\hat{\phi}
\left(\frac{2c_2l\delta}{r^*\sqrt{Q_0}}\right).\label{F2}
\end{equation}
Suppose that $\delta$ satisfies the condition (\ref{P10}). Then,
from (\ref{F1}) and (\ref{F2}), we obtain 
\begin{eqnarray}
& & \frac{1}{\delta} \int\limits_{Q_0}^{2Q_0} \
\sum\limits_{q\in \mathbbm{Z}} \ 
\phi\left(\frac{q-\sqrt{y}}{2c_2\delta/\sqrt{Q_0}} 
\right) \sum\limits_{\scriptsize \begin{array}{cccc} 
m\in \mathbbm{Z}\\ m\equiv -bq^2 \mbox{ mod } r\end{array}} 
\phi\left(\frac{m-yrz}{8\delta rz}\right) {\rm d}y  \label{F3}\\
&\le& \frac{16c_2\delta z}{\sqrt{Q_0}} \sum\limits_{j\in \mathbbm{Z}} 
\frac{\hat{\phi}(8j\delta z)}{r^*} \sum\limits_{l\in \mathbbm{Z}}
\hat{\phi}\left(\frac{2c_2l\delta}{r^*\sqrt{Q_0}}\right)\left\vert 
\sum\limits_{d=1}^{r^*} 
e\left(\frac{j^*bd^2+ld}{r^*}\right) \right\vert \nonumber\\ 
& & \times \left\vert 
\int\limits_{Q_0}^{2Q_0} e\left(jyz-l\cdot\frac{\sqrt{y}}{r^*}\right)\
{\rm d}y\right\vert.\nonumber
\end{eqnarray}
Applying the Lemmas 4 and 7 to the right-hand side of (\ref{F3}), we 
deduce 
\begin{eqnarray}
& & \frac{1}{\delta} \int\limits_{Q_0}^{2Q_0} \
\sum\limits_{q\in \mathbbm{Z}} \ 
\phi\left(\frac{q-\sqrt{y}}{c_2\delta/\sqrt{Q_0}} 
\right) \sum\limits_{\scriptsize \begin{array}{cccc} 
m\in \mathbbm{Z}\\ m\equiv -bq^2 \mbox{ mod } r\end{array}} 
\phi\left(\frac{m-yrz}{8\delta rz}\right) {\rm d}y  \label{F4}\\
&\le& \frac{c_8\delta z}{\sqrt{Q_0}} \sum\limits_{\vert j \vert 
\le 1/(8\delta z)} \frac{1}{\sqrt{r^*}} 
\sum\limits_{\vert l\vert \le r^*\sqrt{Q_0}/(2c_2\delta)} 
\left\vert 
\int\limits_{Q_0}^{2Q_0} e\left(jyz-l\cdot\frac{\sqrt{y}}{r^*}\right)\ {\rm d}y
\right\vert.\nonumber
\end{eqnarray}

If $j=0$ and $l=0$, then the integral on the right-hand side of (\ref{F4}) is 
equal to $Q_0$. If $j\not= 0$ and $l=0$, then 
$$
\left\vert 
\int\limits_{Q_0}^{2Q_0} e\left(jyz-l\cdot\frac{\sqrt{y}}{r^*}\right)\ {\rm d}y
\right\vert \le \frac{1}{\vert j\vert z}.
$$
If $j=0$ and $l\not=0$, then 
$$
\left\vert 
\int\limits_{Q_0}^{2Q_0} e\left(jyz-l\cdot\frac{\sqrt{y}}{r^*}\right)
\ {\rm d}y\right\vert \le \frac{c_{9}Q_0^{1/2}}{\vert l\vert}
$$
by Lemma 5 (take into account that $r^*=1$ if $j=0$).
If $j\not=0$ and $l\not=0$, then Lemma 6 yields
$$
\left\vert 
\int\limits_{Q_0}^{2Q_0} e\left(jyz-l\cdot\frac{\sqrt{y}}{r^*}\right)
\ {\rm d}y\right\vert \le 
\frac{c_{10}\sqrt{r^*}Q_0^{3/4}}{\sqrt{\vert l\vert}}.
$$
Therefore, the right-hand side of (\ref{F4}) can be estimated by
\begin{eqnarray}
&\hspace{-0.5cm}&
 \frac{c_8\delta z}{\sqrt{Q_0}} \sum\limits_{\vert j \vert 
\le 1/(8\delta z)} \frac{1}{\sqrt{r^*}} 
\sum\limits_{\vert l\vert \le r^*\sqrt{Q_0}/(2c_2\delta)} 
\left\vert 
\int\limits_{Q_0}^{2Q_0} e\left(jyz-l\cdot\frac{\sqrt{y}}{r^*}\right)
\ {\rm d}y\right\vert \label{F5}\\
&\le& c_{11}\delta\left(z\sqrt{Q_0}+\frac{1}{\sqrt{Q_0}}
\sum\limits_{1\le j \le 1/(8\delta z)} \frac{1}{j\sqrt{r^*}}\right.
+z\sum\limits_{1\le l\le \sqrt{Q_0}/(2c_2\delta)} \frac{1}{l}+
\nonumber\\ & & \left. zQ_0^{1/4}\sum\limits_{1\le j \le 1/(8\delta z)}\
\sum\limits_{1\le l\le r^*\sqrt{Q_0}/(2c_2\delta)}\frac{1}{\sqrt{l}}\right)
\nonumber\\
&\le& c_{12}\left(\delta z\sqrt{Q_0}+
\frac{\delta}{\sqrt{Q_0}}
\sum\limits_{1\le j \le 1/(8\delta z)} \frac{1}{j\sqrt{r^*}}
+\delta z\Delta^{-\varepsilon}+\right.\nonumber\\
& & \left. z\sqrt{\delta} Q_0^{1/2}
\sum\limits_{1\le j \le 1/(8\delta z)} \sqrt{r^*}\right).\nonumber
\end{eqnarray} 

Now, we evaluate the sums over $j$ in the last line of (\ref{F5}).
By the definition of $r^*$, we have 
\begin{eqnarray}
\sum\limits_{1\le j \le 1/(8\delta z)} \frac{1}{j\sqrt{r^*}}&=&
\frac{1}{\sqrt{r}}
\sum\limits_{t\vert r} \sqrt{t} \sum\limits_{\scriptsize \begin{array}{cccc} 
1\le j \le 1/(8\delta z)\smallskip
\\ (r,j)=t\end{array}} \frac{1}{j}\label{F6} \\
&\le& \frac{c_{13}\log(2+1/(8\delta z))}{\sqrt{r}}
\sum\limits_{t\vert r} \frac{1}{\sqrt{t}}\nonumber\\
&\le& c_{14}\Delta^{-\varepsilon}r^{-1/2}\nonumber
\end{eqnarray}   
and 
\begin{eqnarray}
\sum\limits_{1\le j \le 1/(8\delta z)} \sqrt{r^*}&=&
{\sqrt{r}}
\sum\limits_{t\vert r} \frac{1}{\sqrt{t}} 
\sum\limits_{\scriptsize \begin{array}{cccc} 
1\le j \le 1/(8\delta z)\smallskip\\ (r,j)=t\end{array}} 1 \label{F7} \\
&\le& \frac{\sqrt{r}}{8\delta z}
\sum\limits_{t\vert r} \frac{1}{t^{3/2}}\nonumber\\
&\le& \frac{c_{15}\sqrt{r}}{\delta z}.\nonumber
\end{eqnarray}   

Combining Lemma 2, (\ref{P13}), 
(\ref{F4}), (\ref{F5}), (\ref{F6}) and (\ref{F7}), we
obtain
\begin{equation}
P\left(\frac{b}{r}+z\right)\le 
c_4\Delta^{-\varepsilon}\left(1+\delta z\sqrt{Q_0}
+\delta Q_0^{-1/2}r^{-1/2}+
\delta^{-1/2}Q_0^{1/2}\sqrt{r}\right).\label{Q1}
\end{equation} 
Choosing $\delta:=Q_0\Delta/z$, we infer the desired estimate from 
(\ref{Q1}) and (\ref{P2}).$\Box$

\section{Estimation of $P(b/r+z)$ - second way}
In this section, we use elementary tools to derive the following bound 
for $P(b/r+z)$ from Lemma 2.\\

{\bf Theorem 4:} \begin{it} Suppose that the conditions 
(\ref{P1}), (\ref{P2}) and (\ref{P4}) are satisfied. 
Then,
\begin{equation}
P\left(\frac{b}{r}+z\right)\le c_{16}\Delta^{-\varepsilon}
\left(1+Q_0rz+Q_0^{3/2}\Delta\right).\label{R2}
\end{equation}\end{it}

{\bf Proof:} 
Rearranging the order of summation, the sum on the 
right-hand side of (\ref{PP}) can be written in the form 
\begin{eqnarray}
& & \sum\limits_{\sqrt{y}-c_2\delta/\sqrt{Q_0}\le 
q \le \sqrt{y}+c_2\delta/\sqrt{Q_0}}
\sum\limits_{\scriptsize 
\begin{array}{cccc} (y-4\delta)rz\le m \le (y+4\delta)rz
\smallskip\\ m \equiv -bq^2\mbox{ mod } r\smallskip
\\ m\not=0\end{array}} 1\label{R3}\\
&=& \sum\limits_{\scriptsize \begin{array}{cccc}
(y-4\delta)rz\le m \le (y+4\delta)rz\\ m\not=0\end{array}}
\sum\limits_{\scriptsize \begin{array}{cccc} \sqrt{y}-c_2\delta/\sqrt{Q_0}\le 
q \le \sqrt{y}+c_2\delta/\sqrt{Q_0}\smallskip\\ 
q^2 \equiv -\overline{b}m\mbox{ mod } r\end{array}} 1, \nonumber
\end{eqnarray}
where $\overline{b}$ mod $r$ is the multiplicative inverse of $b$ mod $r$,
{\it i.e.} $\overline{b}b\equiv 1$ mod $r$.
The double sum on the right-hand side of (\ref{R3}) can be split up as 
follows:
\begin{eqnarray}
& & \sum\limits_{\scriptsize \begin{array}{cccc}
(y-4\delta)rz\le m \le (y+4\delta)rz\smallskip\\ m\not=0\end{array}}
\sum\limits_{\scriptsize \begin{array}{cccc} \sqrt{y}-c_2\delta/\sqrt{Q_0}\le 
q \le \sqrt{y}+c_2\delta/\sqrt{Q_0}\smallskip\\ 
q^2 \equiv -\overline{b}m\mbox{ mod } r\end{array}} 1\label{R4}\\
&=& \sum\limits_{t\vert r} 
\sum\limits_{\scriptsize \begin{array}{cccc} 
(y-4\delta)rz/t\le m^{\prime} \le 
(y+4\delta)rz/t\smallskip\\ (m^{\prime},r/t)=1\smallskip\\ 
m^{\prime}\not=0\end{array}}
\sum\limits_{\scriptsize \begin{array}{cccc} 
\sqrt{y}-c_2\delta/\sqrt{Q_0}\le 
q \le \sqrt{y}+c_2\delta/\sqrt{Q_0}\smallskip\\ q^2\equiv 0\mbox{ mod }t
\smallskip\\
q^2/t \equiv -\overline{b}m^{\prime}\mbox{ mod } r/t\end{array}} 1\nonumber\\
\nonumber\\ &=& \sum\limits_{t\vert r} 
\sum\limits_{\scriptsize \begin{array}{cccc} 
(y-4\delta)rz/t\le m^{\prime} \le (y+4\delta)rz/t\smallskip\\ 
(m^{\prime},r/t)=1\smallskip\\ m^{\prime}\not=0\end{array}}
\sum\limits_{\scriptsize \begin{array}{cccc} 
q^{\prime}\in {\cal{S}}_t(y)\smallskip\\ 
q^{\prime} \equiv -\overline{b}m^{\prime}\mbox{ mod } r/t\end{array}} 1,
\nonumber
\end{eqnarray}
where 
$$
{\cal{S}}_t(y):=\left\{q^2/t \ :\  \sqrt{y}-c_2\delta/\sqrt{Q_0}\le 
q \le \sqrt{y}+c_2\delta/\sqrt{Q_0}\ \mbox{ and }\ q^2\equiv 0\mbox{ mod }t
\right\}.
$$
In the following, we determine the structure of ${\cal{S}}_t(y)$. 

Let $t=p_1^{v_1}\cdots p_n^{v_n}$
be the prime number factorization of $t$. For $i=1,...,n$ let
$$
u_i:=\left\{\begin{array}{llll} v_i, & \mbox{ if } v_i \mbox{ is even,}\\ \\
v_i+1, & \mbox{ if } v_i \mbox{ is odd.} \end{array}\right.
$$
Put 
$$
f_t:=p_1^{u_1/2}\cdots p_n^{u_n/2}.
$$
Then $q^2$ $(q\in\mathbbm{N})$ is divisible by $t$ iff $q$ is divisible by
$f_t$. Thus, 
\begin{equation}
{\cal{S}}_t(y)=\left\{q_1^2g_t\ : \ \left(\sqrt{y}-c_2\delta/\sqrt{Q_0}
\right)/f_t \le 
q_1 \le \left(\sqrt{y}+c_2\delta/\sqrt{Q_0}\right)/f_t\right\},\label{R6}
\end{equation}
where  
$$
g_t:=\frac{f_t^2}{t}=p_1^{u_1-v_1}\cdots p_n^{u_n-v_n}.
$$
Hence,
\begin{equation}
\vert {\cal{S}}_t(y)\vert \le 1+\frac{2c_2\delta}{f_t\sqrt{Q_0}}.\label{25}
\end{equation}

By (\ref{R6}),  we get 
\begin{equation}
\sum\limits_{\scriptsize \begin{array}{cccc} 
q^{\prime}\in {\cal{S}}_t(y)\smallskip\\ 
q^{\prime} \equiv -\overline{b}m^{\prime}\mbox{ mod } r/t\end{array}} 1=
\sum\limits_{\scriptsize \begin{array}{cccc} 
\left(\sqrt{y}-c_2\delta/\sqrt{Q_0}
\right)/f_t \le 
q\le \left(\sqrt{y}+c_2\delta/\sqrt{Q_0}\right)/f_t\smallskip\\ 
q^2g_t \equiv -\overline{b}m^{\prime}\mbox{ mod } r/t\end{array}} 1,\label{26}
\end{equation}
and from (\ref{25}) it follows that
\begin{eqnarray}
& & \sum\limits_{\scriptsize \begin{array}{cccc} 
\left(\sqrt{y}-c_2\delta/\sqrt{Q_0}
\right)/f_t\le  
q\le \left(\sqrt{y}+c_2\delta/\sqrt{Q_0}\right)/f_t\smallskip\\ 
q^2g_t \equiv -\overline{b}m^{\prime}\mbox{ mod } r/t\end{array}} 1
\label{27}\\ &\le&  c_{17}s(r,t,m^{\prime})
\left(1+\frac{t\delta}{rf_t\sqrt{Q_0}}\right),\nonumber
\end{eqnarray}
where $s(r,t,m^{\prime})$ is the number of solutions mod 
$r/t$ of the congruence
\begin{equation}
g_tx^2 \equiv -\overline{b}m^{\prime}\mbox{ mod } r/t.\label{29}
\end{equation}

Next, we derive a bound for $s(t,r,m^{\prime})$. In the sequel, we suppose 
that 
\begin{equation}
(m^{\prime},r/t)=1,\label{28}
\end{equation}
as in (\ref{R4}). 
If $(g_t,r/t)>1$, then $s(t,r,m^{\prime})=0$ by (\ref{28}) and 
$(\overline{b},r/t)=1$. Therefore, we can assume that $(g_t,t/r)=1$. Let 
$\overline{g_t}$ mod $r/t$ 
be the multiplicative inverse of
$g_t$ mod $r/t$, {\it i.e.} 
$\overline{g_t}g_t\equiv 1$ mod $r/t$. Put $l:=-\overline{g_tb}m^{\prime}$. 
Then (\ref{29}) is 
equivalent to
\begin{equation}
x^2\equiv l \mbox{ mod } k, \label{30}
\end{equation} 
where $k:=r/t$.
Taking into account that $(l,k)=1$, and using some elementary facts on
the number 
of solutions of polynomial congruences modulo prime powers (see \cite{Qua}, 
for example), we see that  
(\ref{30}) has at most $2$ solutions if $k$ is a power of an odd prime and
at most $4$ solutions if $k$ is a power of 2. From this it follows that
for all $k\in \mathbbm{N}$ and $l\in \mathbbm{Z}$ with $(l,k)=1$ there 
exist 
at most $2^{\omega(k)+1}$ solutions mod $k$ to the congruence (\ref{30}),   
where $\omega(k)$ is the number of distinct 
prime divisors of $k$. Further, we have $2^{\omega(k)}\ll k^{\varepsilon/2}$
(see \cite{Han}). Thus, we obtain
\begin{equation}
s(r,t,m^{\prime})\le c_{18}r^{\varepsilon/2}.\label{31}
\end{equation}

Suppose that the condition (\ref{P10}) is satisfied. Then,
combining (\ref{R3}), (\ref{R4}), (\ref{26}), (\ref{27}) and (\ref{31}), 
we obtain
\begin{eqnarray}
& & \frac{1}{\delta}\int\limits_{Q_0}^{2Q_0}
\sum\limits_{\sqrt{y}-c_2\delta/\sqrt{Q_0}\le 
q \le \sqrt{y}+c_2\delta/\sqrt{Q_0}}
\sum\limits_{\scriptsize 
\begin{array}{cccc} (y-4\delta)rz\le m \le (y+4\delta)rz
\smallskip\\ m \equiv -bq^2\mbox{ mod } r\\m\not=0\end{array}} 1 \ 
{\rm d}y \label{32}\\
&\le& c_{19} r^{\varepsilon/2}
\sum\limits_{t\vert r} \left(1+
\frac{t\delta}{rf_t\sqrt{Q_0}}\right)\cdot\frac{1}{\delta}
\int\limits_{Q_0}^{2Q_0} \sum\limits_{\scriptsize 
\begin{array}{cccc} (y-4\delta)rz/t\le m^{\prime} \le (y+4\delta)rz/t
\smallskip\\ m^{\prime}\not=0\end{array}} 1\ {\rm d}y.
\nonumber
\end{eqnarray}
We estimate the integral on the right-hand side by
\begin{eqnarray}
& & \int\limits_{Q_0}^{2Q_0} \sum\limits_{\scriptsize 
\begin{array}{cccc} (y-4\delta)rz/t\le m^{\prime} \le (y+4\delta)rz/t
\smallskip\\ m^{\prime}\not=0\end{array}} 1\ {\rm d}y\label{L}\\
&=& \sum\limits_{\scriptsize \begin{array}{cccc} 
(Q_0-4\delta)rz/t\le m^{\prime} \le (2Q_0+4\delta)rz/t
\smallskip\\ m^{\prime}\not=0\end{array}} 
\int\limits_{\max\{Q_0,tm^{\prime}/(rz)-4\delta\}}^{\min\{2Q_0,tm^{\prime}/(rz)+4\delta\}} 1
\ {\rm d}y\nonumber\\
&\le& 8\delta \sum\limits_{\scriptsize \begin{array}{cccc} 
(Q_0-4\delta)rz/t\le m^{\prime} \le (2Q_0+4\delta)rz/t
\smallskip\\ m^{\prime}\not=0\end{array}} 1\nonumber\\
&\le & \frac{72\delta Q_0rz}{t}.\nonumber
\end{eqnarray}
To obtain the last line of the above inequality, we use (\ref{P10}) and 
the summation condition $m^{\prime}\not=0$. From (\ref{32}) and (\ref{L}), we
deduce
\begin{eqnarray}
& & \frac{1}{\delta}\int\limits_{Q_0}^{2Q_0}
\sum\limits_{\sqrt{y}-c_2\delta/\sqrt{Q_0}\le 
q \le \sqrt{y}+c_2\delta/\sqrt{Q_0}}
\sum\limits_{\scriptsize 
\begin{array}{cccc} (y-4\delta)rz\le m \le (y+4\delta)rz
\smallskip\\ m \equiv -bq^2\mbox{ mod } r\\m\not=0\end{array}} 1 \ 
{\rm d}y \label{E}\\ &\le& c_{20} r^{\varepsilon/2}
\sum\limits_{t\vert r} \left(1+
\frac{t\delta}{rf_t\sqrt{Q_0}}\right)\frac{Q_0rz}{t}\nonumber\\ &\le&
c_{21} r^{\varepsilon} \left(Q_0rz+\delta z\sqrt{Q_0}\right).
\nonumber
\end{eqnarray}

Finally, we choose 
\begin{equation}
\delta:=\frac{Q_0\Delta}{z},\label{D}
\end{equation}
which is in consistency with the condition (\ref{P10}). 
Combining (\ref{E}), (\ref{D})  and Lemma 2, we obtain the result of 
Theorem 4. $\Box$
 
\section{Proof of Theorem 2}
Suppose that the conditions (\ref{P1}), (\ref{P2}) and (\ref{P4}) are
satisfied. Then the combination of the Theorems 3 and 4 yields

\begin{equation}
P\left(\frac{b}{r}+z\right)\le 
c_{22}\Delta^{-\varepsilon}\left(Q_0^{3/2}\Delta+
\min\left\{Q_0rz,Q_0^{1/2}\Delta r^{-1/2}z^{-1}\right\}+
\Delta^{-1/4}\right).\label{E1}
\end{equation}

If
$$
z\le \Delta^{1/2}Q_0^{-1/4}r^{-3/4},
$$
then 
$$
\min\left\{Q_0rz,Q_0^{1/2}\Delta r^{-1/2}z^{-1}\right\}=
Q_0rz\le Q_0^{3/4}\Delta^{1/2}r^{1/4}.
$$
If 
$$
z> \Delta^{1/2}Q_0^{-1/4}r^{-3/4},
$$
then 
$$
\min\left\{Q_0rz,Q_0^{1/2}\Delta r^{-1/2}z^{-1}\right\}=
Q_0^{1/2}\Delta r^{-1/2}z^{-1}\le Q_0^{3/4}\Delta^{1/2}r^{1/4}.
$$
From the above inequalities and (\ref{P2}), we deduce
\begin{equation}
\min\left\{Q_0rz,Q_0^{1/2}\Delta r^{-1/2}z^{-1}\right\}\le 
Q_0^{3/4}\Delta^{3/8}. \label{E2}
\end{equation}
Furthermore,
\begin{equation}
Q_0^{3/4}\Delta^{3/8}= \sqrt{(Q_0^{3/2}\Delta)\cdot \Delta^{-1/4}}
\le Q_0^{3/2}\Delta+\Delta^{-1/4}. \label{E3}
\end{equation}

Combining (\ref{E1}), (\ref{E2}) and (\ref{E3}), we get
\begin{equation}
P\left(\frac{b}{r}+z\right)\le 
c_{23}\Delta^{-\varepsilon}\left(Q_0^{3/2}\Delta+\Delta^{-1/4}\right).
\label{E4}
\end{equation}
We now choose $\Delta:=1/N$. Then from Theorem 1, Lemma 1 and (\ref{E4}) it 
follows that 
\begin{equation}
\sum\limits_{\sqrt{Q_0}\le q\le \sqrt{2Q_0}} 
\sum\limits_{\scriptsize \begin{array}{cccc} 
a=1\\(a,q)=1\end{array}}^{q^2}
\left\vert S\left(\frac{a}{q^2}\right)\right\vert^2
\ll N^{\varepsilon}\left(Q_0^{3/2}+N^{5/4}\right)Z. 
\label{E5}
\end{equation}
We can devide the interval $[1,Q]$ into $O(\log Q)$ subintervals of the
form $\left[\sqrt{Q_0},\sqrt{2Q_0}\right]$, where $1\le Q_0\le Q^2$. Hence,
the result of Theorem 2 follows from (\ref{E5}). $\Box$
\\ \\ 
 
{\bf Acknowledgement.} This paper was written when the author held a
postdoctoral position at the Harish-Chandra Research Institute at Allahabad 
(India). The author wishes to thank this institute for financial support.\\

\end{document}